\documentclass[11pt]{article}
\usepackage{mathrsfs}
\usepackage{latexsym,lineno}
\usepackage{epsfig}
\usepackage{amsmath}\usepackage{fleqn}\usepackage{verbatim}\usepackage{epsf}
\usepackage{amsthm}\usepackage{graphicx, float}\usepackage{graphicx}
\usepackage{amsfonts}\usepackage{amssymb}\usepackage{graphpap}
\usepackage{epic}\usepackage{curves}

\newtheorem{fact}{Theorem}

\newtheorem{theorem}{Theorem}

\newtheorem{ob}[fact]{Observation}

\newtheorem{lemma}{Lemma}

\title{3-Factor-criticality in double domination edge critical graphs\thanks{This research was partially supported by the National
Nature Science Foundation of China (No. 11171207), the
Nature Science Foundation of Shanghai (No. 14ZR1417900) and the Key Programs of
Wuxi City College of Vocational Technology(WXCY-2012-GZ-007).} }
\author{Haichao Wang$^{1}$,  \, Erfang Shan$^{2}$\thanks {
Corresponding author. Email address: efshan@shu.edu.cn (E. Shan)}, \, Yancai Zhao$^{3}$\\
{\small $^{1}$School of Mathematics and Physics, Shanghai
University of Electric Power,}\\
{\small Shanghai 200090, China}\\
{\small $^{2}$School of Management, Shanghai University, Shanghai 200444, China}\\
{\small $^{3}$Department of Basic Science, Wuxi City College of Vocational Technology,}\\
{\small Jiangsu 214153, China}\\
 }

\date{}
\begin{document}
\maketitle \baselineskip17pt
\begin{abstract}
\baselineskip17pt

 A vertex subset $S$ of a graph $G$ is a double dominating set of
 $G$ if $|N[v]\cap S|\geq 2$ for each vertex $v$ of $G$, where
 $N[v]$ is the set of the vertex $v$  and vertices adjacent to $v$.  The double domination
 number of $G$, denoted by $\gamma_{\times 2}(G)$, is the cardinality of a
 smallest double dominating set of $G$. A graph $G$
 is said to be double domination edge critical if $\gamma_{\times
 2}(G+e)<\gamma_{\times 2}(G)$ for any edge $e \notin E$. A double
 domination edge critical graph $G$ with $\gamma_{\times 2}(G)=k$ is
 called $k$-$\gamma_{\times 2}(G)$-critical. A graph $G$ is $r$-factor-critical
 if $G-S$ has a perfect matching for each set $S$ of $r$ vertices in $G$. In this paper we show
 that $G$ is 3-factor-critical if $G$ is a 3-connected claw-free $4$-$\gamma_{\times
 2}(G)$-critical graph of odd order with minimum degree at least 4 except a family of graphs.

\bigskip


\bigskip

\end{abstract}

\section{Introduction}

Recently, the matching and factor properties in critical graphs
with respect to domination
 have received more attention (see, [1-5, 7, 10, 14-15, 21-24]).  A graph $G$ is {\em $r$-factor-critical} if $G-S$
has a perfect matching for each set $S$ of $r$ vertices in $G$.
If $r=1$, the graph is said to be {\em factor-critical} and if
$r=2$, the graph $G$ is called {\em bicritical}. A
 {\em double dominating set} (DDS) of $G$ is defined in \cite{hh} as a subset $S$ of $V(G)$
 such that $|N[v]\cap S|\geq 2$ for every vertex $v$ of $G$, where $N[v]$
 is the set of the vertex $v$ and vertices adjacent to $v$ in $G$. The
 {\em double domination number} $\gamma_{\times 2}(G)$ of $G$
 is the cardinality of a smallest double dominating set of $G$.  A graph $G$ is called {\em double
domination edge critical}, or just {\em $\gamma_{\times
2}(G)$-critical}, if $\gamma_{\times
 2}(G+e)<\gamma_{\times 2}(G)$ for each edge $e \not\in E(G)$.
If $\gamma_{\times 2}(G)=k$, a $\gamma_{\times 2}(G)$-critical graph
 is said to be $k$-$\gamma_{\times 2}(G)$-{\em
critical}.  In \cite{wk,ws}  the matching properties of double domination edge
critical graphs were investigated, we proved that $G$ has a perfect matching if $G$ is a connected $K_{1,4}$-free 4-$\gamma_{\times 2}(G)$-critical
graph of even order $\ge 6$ with an exceptional family of graphs;  $G$ is bicritical if
$G$ is either a 2-connected claw-free 4-$\gamma_{\times 2}(G)$-critical of even order with minimum degree at
least 3 or  a 3-connected $K_{1,4}$-free 4-$\gamma_{\times 2}(G)$-critical
graph of even order with minimum degree at least 4.
In this paper we show
 that $G$ is 3-factor-critical if $G$ is a 3-connected claw-free $4$-$\gamma_{\times
 2}(G)$-critical graph of odd order with minimum degree at least 4 except a family of graphs.

 For notation and graph theory terminology we in general follow
 \cite{bm}. Specifically, let $G$ be a finite simple graph with
{\em vertex set} $V(G)$ and {\em edge set} $E(G)$.
 For a vertex $v\in V(G)$, the {\em open neighborhood} of
$v$ is $N(v)=\{u\in V(G)\colon\, uv\in E(G)\}$ and the {\em closed
neighborhood} of $v$ is $N[v]=\{v\}\cup N(v)$. The {\em degree} of
$v$ in $G$, denoted by $d(v)$, is the cardinality of $N(v)$. Let
$\delta(G)$ represent the {\em minimum degree} of $G$. As usual,  $K_{m,n}$ denotes a complete bipartite graph with  classes of
cardinality $m$ and $n$; $K_n$  is the  complete graph on $n$ vertices, and  $C_n$
is the cycle on $n$ vertices.
For
$S\subseteq V(G)$, the subgraph of $G$ induced by $S$ is denoted by
$G[S]$. A graph $G$ is said to be {\em $K_{1,r}$-free} if it
contains no  $K_{1,r}$ as an induced
subgraph. In particular, $K_{1,3}$-free is also called {\em
claw-free}. The {\em complement} of $G$, denoted by
$\overline{G}$, is the graph with vertex set $V(G)$ such that two
vertices are adjacent in $\overline{G}$ if and only if the
vertices are not adjacent in $G$. The {\em diameter} of $G$ is the greatest distance between
two vertices of $G$, denoted by diam$(G)$.
A {\em cutset} of a connected $G$ is a subset $S$ of $V(G)$ such that $G-S$
is disconnected. For $S\subseteq
V(G)$, we shall denote by $\omega(G-S)$, the number of components of
$G-S$ and by $o(G-S)$, the number of odd components of $G-S$.
A subset $S$ of $V(G)$ is called an {\em independent set} of $G$ if
no two vertices of $S$ are adjacent in $G$. The {\em independence
number} of $G$, denoted by $\alpha(G)$, is the cardinality of a
largest independent set of $G$. A set of pairwise independent
edges in  $G$ is called a {\em matching} of $G$. A matching is
{\em perfect }if it is incident with every vertex of $G$.

 For a fixed positive integer $k$, a
 $k$-{\em tuple dominating set} of $G$ is a subset $S$ of $V(G)$
 such that $|N[v]\cap S|\geq k$ for every vertex $v\in V(G)$. The
 $k$-{\em tuple domination number} $\gamma_{\times k}(G)$ of $G$
 is the minimum cardinality of a $k$-tuple dominating set of $G$. In particular,
 when $k=1,2$,  1-tuple domination and 2-tuple domination are the
 ordinary domination and double domination, respectively. The concept of $k$-tuple domination in graphs was
 introduced and studied in \cite{hh}. For more
results on the $k$-tuple domination, we refer to [6, 9, 11-13,
16-18, 19-21, 25].

\section{Preliminaries}

In this section we  state some  results that are useful in the
proof of our main results. For an edge $uv\in E(\overline{G})$, we
shall denote by $D_{uv}$ a minimum double dominating set (MDDS) of $G+uv$ throughout this
paper.

By the definition of $\gamma_{\times2}$-criticality, the
following observation follows immediately.

\begin{ob}\label{Observation1}
If $G$ is a $\gamma_{\times2}$-critical graph and $uv\in
E(\overline{G})$, then $D_{uv}$ contains at least one of $u$ and
$v$. Furthermore, if $\gamma_{\times 2}(G+uv)=\gamma_{\times
2}(G)-2$, then $D_{uv}$ contains both $u$ and $v$.
\end{ob}

\begin{lemma}{\rm (\cite{t})}\label{lem1}
If $G$ is a connected $4$-$\gamma_{\times 2}(G)$-critical graph,
then diam$(G)\in\{2,3\}$.
\end{lemma}

\begin{lemma}{\rm (\cite{t})}\label{lem2}
A graph $G$ with diam$(G)=3$ is $4$-$\gamma_{\times
2}(G)$-critical if and only if $G$ is the sequential join
$K_{1}+K_{s}+K_{t}+K_{1}$ for positive integers $s$ and $t$.
\end{lemma}

The {\em sequential join}, as defined by Akiyama and Harary, for
three or more disjoint graphs $G_{1},G_{2},\ldots, G_{n}$, denoted
by $G_{1}+G_{2}+\cdots+ G_{n}$, is the graph
$(G_{1}+G_{2})\cup(G_{2}+G_{3})\cup\cdots\cup(G_{n-1}+ G_{n})$,
where $G_{i}+G_{i+1}$ is obtained from $G_i\cup G_{i+1}$ by
joining each vertex of $G_i$ to each vertex of $G_{i+1}$ for $1\le
i\le n-1$.

In \cite{ws} Wang and Shan  proved the following two results.

\begin{lemma}{\rm (\cite{ws})}\label{lem3}
If $G$ is a connected $K_{1,r}$-free {\em($k\geq3$)}
$4$-$\gamma_{\times 2}(G)$-critical graph, then $\alpha(G)\leq r$.
\end{lemma}

\begin{lemma}{\rm (\cite{ws})}\label{lem4}
Let $G$ be a connected $4$-$\gamma_{\times2}(G)$-critical graph and $S$  a cutset
of $G$. If $\omega(G-S)\geq 3$ and $x$ and $y$ belong to different
components of $G-S$, then $|D_{xy}\cap\{x,y\}|=1$ and
$|D_{xy}|=3$.
\end{lemma}

By  Observation \ref{Observation1}, we immediately have the following lemma.

\begin{lemma}\label{lem5}
Let $G$ be a connected $4$-$\gamma_{\times2}(G)$-critical graph and $S$  a cutset
of $G$. If $\omega(G-S)=2$ and each component has at least two
vertices, and $x$ and $y$ lie in different components of $G-S$,
then $|D_{xy}|=3$ and $D_{xy}\cap S\neq \emptyset $.
\end{lemma}

Furthermore, the following results are useful in the proof of our main
result. The first result is proved by Wang and Kang in \cite{wk} and the
second result is due to Favaron \cite{f}.

\begin{theorem}{\rm (\cite{wk})}\label{the0}
If $G$ is a connected
 $K_{1,4}$-free $4$-$\gamma_{\times2}(G)$-critical graph of odd
 order and $\delta(G)\geq 2$, then $G$ is factor-critical.
\end{theorem}

\begin{theorem}{\rm (\cite{f})}\label{the00}
A graph $G$ is $k$-factor-critical if and only if $o(G-S)\leq
|S|-k$ for every $S\subseteq V(G)$ and $|S|\geq k$.
\end{theorem}

\section{Main result}

In this section we shall show that $G$ is 3-factor-critical if $G$
is a 3-connected claw-free $4$-$\gamma_{\times
 2}(G)$-critical graph of odd order with minimum degree at least 4
except a family $\mathcal H$ of graphs.

For convenience, let us introduce more notation and terminology. If $S\subseteq V(G)$ is a minimum double dominating
 set (DDS) of $G$, we call $S$ a $\gamma_{\times 2}(G)$-set. For a vertex
 $v\in V(G)$, $v$ is said to be {\em doubly dominated} by $S$ if
 $|N[v]\cap S|\geq2$. For $A,B\subseteq V(G)$, we say that $A$ is
 {\em doubly dominated} by $B$, written $B\succ_{\times2} A $,
  if each vertex of $A$ is doubly dominated by $B$. Furthermore, we use
 $B\nsucc_{\times2} A$ to present that $A$ is not doubly dominated by $B$.

\begin{lemma}\label{lem6}
Let $G$ be a $4$-$\gamma_{\times
2}(G)$-critical graph of odd order with $\delta(G)\geq 4$.
If  diam$(G)=3$, then $G$ is $3$-factor-critical.
\end{lemma}
\noindent{\bf Proof.}  If diam$(G)=3$,  then, by Lemma \ref{lem2}, $G$ is isomorphic to a
sequential join $K_{1}+K_{s}+K_{t}+K_{1}$ for positive integers
$s$ and $t$. Since $\delta(G)\geq 4$ and $G$ has odd order, it
follows that $s\geq 4$ and $t\geq 4$. Further, $s$ and $t$ must
have different parities. It is easy to check that $G-D$ has a
perfect matching for each set $D$ of $3$ vertices in $G$. So the
assertion holds. ~$\Box$

\begin{lemma}\label{lem7}
Let $G$ be a $3$-connected claw-free $4$-$\gamma_{\times
2}(G)$-critical graph of odd order with $\delta(G)\geq 4$.
If $G$ is not $3$-factor-critical, then there exists a subset $S\subseteq V(G)$
such that $|S|=3$ and $G-S$ contains exactly two odd components.
\end{lemma}
\noindent{\bf Proof.} Since $G$ is not 3-factor-critical,
there exists a subset $S\subseteq V(G)$ with $|S|\geq 3$ such that
$o(G-S)>|S|-3$ by Theorem \ref{the00}. But, by Theorem \ref{the0},
$G$ is factor-critical, and so $o(G-S)\leq |S|-1$. Note that
$|V(G)|$ is odd, so $o(G-S)=|S|-1$ by parity. Thus  $|S|\leq 4$ by Lemma \ref{lem3}.
 Since $G$ is
3-connected, $3\leq |S|\leq 4$.

If $|S|=4$, by Lemma \ref{lem3}, then $G-S$ has no even components
and $\omega(G-S)=3$. Let $S=\{u_{1},u_{2},u_{3},u_{4}\}$. Choose
$a\in V(C_{1})$ and $b\in V(C_{2})$. Now consider the graph
$G+ab$. By Lemma \ref{lem4}, $|D_{ab}\cap\{a,b\}|=1$ and
$|D_{ab}|=3$. Without loss of generality, assume that $a\in
D_{ab}$. In order to doubly dominate $V(C_2)\cup V(C_{3})$, $1\le
|D_{ab}\cap S|\le 2$. If $|D_{ab}\cap S|=2$, without loss of
generality, say $u_{1},u_{2}\in D_{ab}$, then each vertex of
$V(C_{2})\cup V(C_{3})-\{b\}$ is adjacent to both $u_{1}$ and $u_{2}$.
Moreover, $b$ is adjacent to only one of $u_{1}$ and $u_{2}$ for
otherwise $\{a,u_{1},u_{2}\}$ would be a DDS of
$G$, a contradiction. Then $|V(C_{2})|\geq 3$ as
$d(b)\ge\delta(G)\geq 4$. Since $D_{ab}=\{a,u_{1},u_{2}\}$ is a
$\gamma_{\times 2}(G+ab)$-set, $a$ is adjacent to at least one of
$u_{1}$ and $u_{2}$. This implies that $G$ contains a claw
centered at $u_{1}$ or $u_{2}$, a contradiction. Hence
$|D_{ab}\cap S|=1$. This implies that $|D_{ab}\cap V(C_{3})|=1$
and $V(C_2)=\{b\}$. Without loss of generality, let $u_{1}\in
D_{ab}$ and $c\in D_{ab}$ where $c\in V(C_{3})$. Then
$G[\{u_{1},a,b,c\}]$ is a claw in $G$, a contradiction.
Therefore, $|S|\leq 3$. Then we have
$|S|=3$. Moreover, since $S$ is a minimum cutset of $G$,
it follows that each vertex of $S$ is adjacent to
a vertex of each component of $G-S$. Recall that $G$ is claw-free.
Thus $G-S$ has no even components, and so $G-S$ contains exactly two odd components. ~$\Box$

Let $G$ be defined as that in Lemma \ref{lem7}, $S=\{s_{1},s_{2},s_{3}\}$ and
$\{C_{1}, C_2\}$ be two odd components of $G-S$.
Since $\delta(G)\geq 4$, $|V(C_{1})|\geq 3$ and $|V(C_{2})|\geq 3$.  Now set
$A_{i}= V(C_{1})\cap N(s_{i})$ and $B_{i}= V(C_{2})\cap N(s_{i})$ for
$1\leq i\leq3$. We have the following lemma.

\begin{lemma}\label{lem8}
Let $G$ be a $3$-connected claw-free $4$-$\gamma_{\times
2}(G)$-critical graph of odd order with $\delta(G)\geq 4$.
If $G$ is not $3$-factor-critical and diam$(G)=2$, then the following statements are true:\\
$(1) $  For $1\leq i\leq3$, $A_{i}\neq
\emptyset$ and $B_{i}\neq \emptyset$. Furthermore, both $G[A_{i}]$ and $G[B_{i}]$
are complete;\\
$(2)$ $V(C_{1})=\cup_{i=1}^{3}A_{i}$ and $V(C_{2})=\cup_{i=1}^{3}B_{i}$;\\
$(3)$ There exists at least a pair of $A_{i}$ and $ A_{j}$ such that $A_{i}\cap A_{j}\neq\emptyset$, where $1\leq
i\neq j\leq 3$;\\
$(4)$ There exists at least a pair of $B_{i}$ and $ B_{j}$ such that $B_{i}\cap B_{j}\neq\emptyset$, where $1\leq
i\neq j\leq 3$.\\
\end{lemma}
\noindent{\bf Proof.} (1) The statement (1) directly follows, because $G$ is claw-free and $S$ is a minimum cutset of $G$.

(2) Suppose not, without loss of generality, let
$V(C_{1})\neq\cup_{i=1}^{3}A_{i}$. Thus there exists a vertex $u\in V(C_{1})-\cup_{i=1}^{3}A_{i}$. Take any
vertex $v\in V(C_{2})$. Clearly, the distance between $u$ and $v$ is
more than 2, which contradicts our assumption that diam$(G)=2$. So
the statement (2) holds. 

(3) Suppose to the contrary that $A_{i}\cap A_{j}=\emptyset$ for $1\leq i\neq j\leq 3$. Choose
$a_{1}\in A_{1}$ and $b_{1}\in B_{1}$. By Observation 1, $1\leq
|D_{a_{1}b_{1}}\cap\{a_{1},b_{1}\}|\leq2 $. If
$|D_{a_{1}b_{1}}\cap\{a_{1},b_{1}\}|=2$, by Lemma \ref{lem5}, then
$D_{a_{1}b_{1}}\cap S=\{s_{1}\}$ as $A_{i}\cap A_{j}=\emptyset$
for $1\leq i\neq j\leq 3$. However, $A_{2}$ and $A_{3}$ are not doubly dominated by $D_{a_{1}b_{1}}$,
 a contradiction. Hence
$|D_{a_{1}b_{1}}\cap\{a_{1},b_{1}\}|=1$. If $b_{1}\in D_{a_{1}b_{1}}$, then $s_{1} \in D_{a_{1}b_{1}}$,
so that $a_{1}$ can be doubly dominated. But this implies that $A_{2}$ and $A_{3}$ can not be doubly dominated by $D_{a_{1}b_{1}}$,
a contradiction.
Thus $a_{1}\in D_{a_{1}b_{1}}$. Clearly,
$1\leq |D_{a_{1}b_{1}}\cap S|\leq2$. If $|D_{a_{1}b_{1}}\cap
S|=1$, then $s_{1}\in D_{a_{1}b_{1}}$ because $A_{i}\cap
A_{j}=\emptyset$ for $1\leq i\neq j\leq 3$. To doubly dominate
$V(C_{2})-\{b_{1}\}$, we have $|D_{a_{1}b_{1}}\cap(V(C_{2})-\{b_{1}\})|=1$.
Let $x\in D_{a_{1}b_{1}}\cap(V(C_{2})-\{b_{1}\})$. Then $s_{1}x\in E(G)$. By Claim 1,
$D_{a_{1}b_{1}}$ is also a DDS of $G$, a contradiction. Hence
$|D_{a_{1}b_{1}}\cap S|=2$. Since $a_{1}s_{2},a_{1}s_{3}\notin
E(G)$, $s_{1}\in D_{a_{1}b_{1}}$. Without loss of generality,
assume that $s_{2}\in D_{a_{1}b_{1}}$. For a vertex $a_{3}\in
A_{3}$, then $a_{3}$ is adjacent to at least one of $s_{1}$ and
$s_{2}$ to doubly dominate $a_{3}$. This implies that $A_{1}\cap
A_{3}\neq\emptyset$ or $A_{2}\cap A_{3}\neq\emptyset$, a
contradiction. So the statement (3) follows.

(4) We can show that the statement is also true by a similar argument that used in the proof of the statement (3). ~$\Box$

The family $\mathcal H$ of
graphs is defined as follows:  For odd integer $r\geq 3$, let
$H_{1}=K_{r}$, $H_{2}=K_{3}$ and $H_{3}=K_1\cup K_2$. Let $H_{r,3,3}$ be the graph obtained
from $(H_{1}+H_{3})\cup H_{2}$ by adding 6 edges between
$V(H_{2})$ and $V(H_{3})$ such that each vertex of $H_{2}$ has
exactly 2 neighbors in $H_{3}$ while each vertex of $H_{3}$ has
precisely 2 neighbors in $H_{2}$. By our construction, it is easy
to verified that $H_{r,3,3}$ is a 3-connected claw-free
$4$-$\gamma_{\times2}$-critical graphs of odd order with minimum
degree 4. Obviously, $H_{r,3,3}-V(H_{3})$ has no perfect matching, hence $H_{r,3,3}$ is not 3-factor-critical. Let ${\mathcal
H}=\{H_{r,3,3}\colon\, r\geq 3\ \mbox{is an odd integer}\}$.

\begin{lemma}\label{lem9}
Let $G$ be a $3$-connected claw-free $4$-$\gamma_{\times
2}(G)$-critical graph of odd order with $\delta(G)\geq 4$ and diam$(G)=2$.
If $G$ is not $3$-factor-critical and $G\notin{\mathcal H}$, then
both $\cap_{i=1}^{3}A_{i}=\emptyset$ and $\cap_{i=1}^{3}B_{i}=\emptyset$.
\end{lemma}
\noindent{\bf Proof.} Let $S$, $C_{1}$ and $C_{2}$ be defined as before.
 Suppose not, without loss of generality, let
$\cap_{i=1}^{3}A_{i}\neq\emptyset$.
Take $a\in \cap_{i=1}^{3}A_{i}$. By Lemma \ref{lem8} (1), $N[a]=S\cup
V(C_{1})$. We claim that $\cap_{i=1}^{3}B_{i}=\emptyset$.
Otherwise, there exists a vertex $b\in \cap_{i=1}^{3}B_{i}$.
Consider $G+ab$. By Observation 1, without loss of generality,
assume that $a\in D_{ab}$. If $b\in D_{ab}$,
then $|D_{ab}\cap S|=1$ by Lemma \ref{lem5}, so $D_{ab}$ is also a DDS of $G$
because $a $ and $b$ are adjacent to every vertex of $S$, a
contradiction. Thus $b\notin D_{ab}$. Since $D_{ab}\succ_{\times2}
S\cup( V(C_{2})-\{b\}) $, $|D_{ab}\cap (S\cup(
V(C_{2})-\{b\}))|=2$. Since $b\in \cap_{i=1}^{3}B_{i}$, $D_{ab}$
is still a DDS of $G$ by Lemma \ref{lem8} (1), a contradiction. Therefore,
$\cap_{i=1}^{3}B_{i}=\emptyset$. By Lemma \ref{lem8} (4), we may assume that
$B_{1}\cap B_{2}\neq\emptyset$. Let $x\in B_{1}\cap B_{2}$ and
consider the graph $G+ax$.

\noindent{\bf Case 1.} $|D_{ax}\cap \{a,x\}|=2$.

By Lemma \ref{lem5}, we have $s_{3}\in D_{ax}$, and so $
D_{ax}=\{a,x,s_{3}\}$. Then each vertex of $V(C_{1})\cup
(V(C_{2})-\{x\})$ is adjacent to $s_{3}$. Choose $b_{1}\in
V(C_{2})-\{x\}$. Now consider $G+ab_{1}$.

\noindent{\bf Case 1.1.} $|D_{ab_{1}}\cap \{a,b_{1}\}|=2$.

Without loss of generality, suppose that $D_{ab_{1}}= \{a,s_{1},b_{1}\}$. Then each vertex of
$V(C_{1})$ is adjacent to $s_{1}$ while each vertex of $V(C_{2})-\{b_{1}\}$ is adjacent to both $s_{1}$ and
$b_{1}$. Choose a vertex $u\in V(C_{2})-\{x,b_{1}\}$. This derives that $\{s_{1},s_{3},u\}$ is a DDS of $G$, a
contradiction. So $|D_{ab_{1}}\cap \{a,y_{1}\}|=2$ is impossible.

\noindent{\bf Case 1.2.} $|D_{ab_{1}}\cap \{a,b_{1}\}|=1$.

Since $N[a]=V(C_{1})\cup S$, we have $a\in D_{ab_{1}}$ and $b_{1}\notin D_{ab_{1}}$.
Assume that $|D_{ab_{1}}\cap S|=2$. Then $D_{ab_{1}}=\{a,s_{1},s_{3}\}$ or $\{a,s_{2},s_{3}\}$ because
$\cap_{i=1}^{3}B_{i}=\emptyset$. However, $D_{ab_{1}}\nsucc_{\times2} \{x\}$ as $xs_{3}\notin E(G)$, a contradiction.
Hence $|D_{ab_{1}}\cap S|=1$. Notice that $xs_{3}\notin E(G)$. We deduce that $s_{1}\in D_{ab_{1}}$ or $s_{2}\in D_{ab_{1}}$.
Without loss of generality, suppose that $s_{1}\in D_{ab_{1}}$. To doubly dominate $V(C_{2})-\{b_{1}\}$, we see that
$|D_{ab_{1}}\cap (V(C_{2})-\{b_{1}\})|=1$. Let $ b_{2}\in D_{ab_{1}}\cap (V(C_{2})-\{b_{1}\})$. Thus each vertex of
$V(C_{1})$ is adjacent to $s_{1}$ while each vertex of $V(C_{2})-\{b_{1}\}$ is adjacent to both $s_{1}$ and
$b_{2}$. If $b_{1}s_{1}\in E(G)$, then $\{s_{1},s_{3},b_{2}\}$ is a DDS of $G$, a contradiction. So $b_{1}s_{1}\notin E(G)$.
Thus, in order to doubly dominate $b_{1}$ in $G+ab_{1}$,
$b_{1}b_{2}\in E(G)$. This implies that
$\{s_{1},s_{3},b_{2}\}$ is still a DDS of $G$, a contradiction again.

\noindent{\bf Case 2.} $|D_{ax}\cap \{a,x\}|=1$.

If $x\in D_{ax}$, then $D_{ax}$ is also a DDS of $G$ as
$N[a]=S\cup V(C_{1})$, a contradiction. Hence $a\in D_{ax}$. By Lemma \ref{lem5},
$1\leq|D_{ax}\cap S|\leq2$. Suppose that $|D_{ax}\cap S|=2$. Then
$D_{ax}=\{a,s_{1},s_{3}\}$ or $\{a,s_{2},s_{3}\}$. If
$D_{ax}=\{a,s_{1},s_{3}\}$, then each vertex of $V(C_{2})-\{x\}$
is adjacent to both $s_{1}$ and $s_{3}$, i.e.,
$(V(C_{2})-\{x\})\subseteq B_{1}\cap B_{3}$. Since
$\cap_{i=1}^{3}B_{i}=\emptyset$, $s_{2}$ is not adjacent to any
vertex of $V(C_{2})-\{x\}$. Choose $x_{1}\in V(C_{2})-\{x\}$. Now
consider $G+ax_{1}$. If $|D_{ax_{1}}\cap \{a, x_{1}\}|=2$, then
$D_{ax_{1}}=\{a,s_{2},x_{1}\}$ by Lemma \ref{lem5}, however,
$D_{ax_{1}}\nsucc_{\times2}V(C_{2})-\{x,x_{1}\}$. So
$|D_{ax_{1}}\cap \{a, x_{1}\}|=1$. Since $N[a]=S\cup V(C_{1})$,
$x_{1}\notin D_{ax_{1}}$. Then $a\in D_{ax_{1}}$. To doubly
dominate $V(C_{2})-\{x,x_{1}\}$, it follows that $s_{2}\notin
D_{ax_{1}}$ and $|D_{ax_{1}}\cap \{s_{1},s_{3}\}\cup
(V(C_{2})-x_{1})|=2$. Then, since $V(C_{2})=B_{1}$ and $G[B_{1}]$ is complete,
 $D_{ax_{1}}$ is a DDS of $G$, this is a contradiction. If $D_{ax}=\{a,s_{2},s_{3}\}$,
 then we can reach a contradiction by
similar arguments. This implies that
$|D_{ax}\cap S|=2$ is impossible. Hence $|D_{ax}\cap S|=1$.
Suppose $D_{ax}$ contains $s_{1}$ or $s_{2}$. Then $|D_{ax}\cap
(V(C_{2})-\{x\})|=1$, so that $V(C_{2})-\{x\}$ can be doubly dominated. Let
$x_{2}\in D_{ax}\cap (V(C_{2})-\{x\})$. Clearly, $xx_{2}\notin E(G)$
for otherwise $D_{ax}$ would be a DDS of $G$. But then a claw would
occur at $s_{1}$ or $s_{2}$ in $G$, a contradiction. Hence,
$s_{3}\in D_{ax}$. To doubly dominate $V(C_{2})-\{x\}$, we have
$|D_{ax}\cap (V(C_{2})-\{x\})|=1$. Without loss of generality, let
$x_{3}\in D_{ax}\cap (V(C_{2})-\{x\})$. Then $
D_{ax}=\{a,s_{3},x_{3}\}$. Thus $xx_{3}\in E(G)$ and each vertex
of $V(C_{1})$ is adjacent to $s_{3}$ while each vertex of
$V(C_{2})-\{x,x_{3}\}$ is adjacent to both $x_{3}$ and $s_{3}$.

Now we consider $G+ax_{3}$. By Lemma \ref{lem5}, $|D_{ax_{3}}|=3$ and
$D_{ax_{3}}\cap S\neq\emptyset$.

\noindent{\bf Case 2.1.} $|D_{ax_{3}}\cap \{a,x_{3}\}|=2$.

In this subcase, we have $|D_{ax_{3}}\cap \{s_{1},s_{2}\}|=1$. Without loss of
generality, assume that $s_{2}\in D_{ax_{3}}$. Then
$D_{ax_{3}}=\{a,s_{2},x_{3}\}$ and each vertex of
$V(C_{1})$ is adjacent to $s_{2}$ while each vertex of $V(C_{2})-\{x_{3}\}$ is adjacent to both $s_{2}$ and $x_{3}$.
Furthermore, $s_{2}x_{3}\notin E(G)$.
Note that $s_{2}s_{3}\in E(G)$, because $D_{ax}=\{a,s_{3},x_{3}\}$ and $D_{ax}$
doubly
dominates $s_{2}$ in $G+ax$. Choose $y\in V(C_{2})-\{x,x_{3}\}$ and
consider $G+ay$. Suppose $|D_{ay}\cap\{a,y\}|=1$. Then $a\in
D_{ay}$ because $N[a]=S\cup V(C_{1})$. If $|D_{ay}\cap S|=2$, then
$D_{ay}=\{a,s_{1},s_{2}\}$ or $\{a,s_{1},s_{3}\}$. If
$D_{ay}=\{a,s_{1},s_{2}\}$, then $D_{ay}\nsucc_{\times2}\{x_{3}\}$
as $s_{2}x_{3}\notin E(G)$. So $D_{ay}=\{a,s_{1},s_{3}\}$. Since
$\cap_{i=1}^{3}B_{i}=\emptyset$, $xs_{3}\notin E(G)$, and so
$D_{ay}\nsucc_{\times2}\{x\}$. Hence $|D_{ay}\cap S|=1$. By
similar arguments above, we obtain $D_{ay}\cap
\{s_{2},s_{3}\}=\emptyset$. Hence $s_{1}\in D_{ay}$. Since
$\cap_{i=1}^{3}B_{i}=\emptyset$, if $|V(C_{2})|\geq 5$, then
$D_{ay}\nsucc_{\times2}V(C_{2})-\{x,y,x_{3}\}$, a contradiction.
So $|V(C_{2})|=3$, i.e., $V(C_{2})=\{x,y,x_{3}\}$. To doubly
dominate $y$, we have $|D_{ay}\cap\{x,x_{3}\}|=1$. If $x\in D_{ay}$, then
$x_{3}s_{1}, s_{1}s_{3}\in E(G)$, which implies that $S$ is a DDS
of $G$, a contradiction. Hence $x_{3}\in D_{ay}$. To doubly
dominate $s_{2}$,  we see that $s_{1}s_{2}\in E(G)$ as $s_{2}x_{3}\notin E(G)$.
This means that $S$ is a DDS of $G$, a contradiction. Therefore,
$|D_{ay}\cap\{a,y\}|=2$. Then $D_{ay}\cap S=\{s_{1}\}$ by Lemma
\ref{lem5}. Thus each vertex of $V(C_{1})$ is adjacent to $s_{1}$
while each vertex of $V(C_{2})-\{y\}$ is adjacent to both $y$ and
$s_{1}$. If $|V(C_{2})|\geq 5$, then
$\cap_{i=1}^{3}B_{i}\neq\emptyset$, a contradiction. So
$|V(C_{2})|=3$, i.e., $V(C_{2})=\{x,y,x_{3}\}$. Clearly,
$ys_{1}\notin E(G)$. Further, we can obtain
$s_{1}s_{2},s_{1}s_{3}\notin E(G)$ for otherwise $S$ would be a
DDS of $G$. By Lemma \ref{lem8} (1), $G$ is isomorphic to $H_{r,3,3}$, a
contradiction. Hence $|D_{ax_{3}}\cap \{a,x_{3}\}|=2$ is
impossible.

\noindent{\bf Case 2.2.} $|D_{ax_{3}}\cap \{a,x_{3}\}|=1$.

In this subcase, clearly  $a\in D_{ax_{3}}$ and $x_{3}\notin
D_{ax_{3}}$. Since $xs_{3},ax\notin E(G)$, $s_{3}\notin
D_{ax_{3}}$. If $|D_{ax_{3}}\cap S|=2$, then
$D_{ax_{3}}=\{a,s_{1},s_{2}\}$. Thus each vertex of
$V(C_{2})-\{x,x_{3}\}$ is adjacent to both $s_{1}$ and $s_{2}$, and so
$\cap_{i=1}^{3}B_{i}=V(C_{2})-\{x,x_{3}\}\neq\emptyset$, a
contradiction. Hence $|D_{ax_{3}}\cap S|=1$. Without loss of
generality, suppose that $s_{2}\in D_{ax_{3}}$. Then each vertex
of $V(C_{1})\cup( V(C_{2})-\{x_{3}\})$ is adjacent to $s_{2}$ and
$( V(C_{2})-\{x,x_{3}\})\subseteq B_{2}\cap B_{3}$. Since $x_{3}$
is adjacent to each vertex of $V(C_{2})-\{x_{3}\}$,
$s_{2}x_{3}\notin E(G)$ for otherwise $D_{ax_{3}}$ is a DDS of
$G$. Thus $s_{2}s_{3}\in E(G)$ so that $s_{2}$ can be doubly dominated  in
$G+ax$. Choose $y_{1}\in V(C_{2})-\{x,x_{3}\}$ and consider
$G+ay_{1}$. By an argument similar to that as in Case
2.1, one can arrive at a contradiction. Therefore,
$|D_{ax_{3}}\cap \{a,x_{3}\}|=1$ is also impossible. ~$\Box$

\begin{theorem}\label{the1}
Let $G$ be a $3$-connected claw-free $4$-$\gamma_{\times
2}(G)$-critical graph of odd order with $\delta(G)\geq 4$.
If  $G\notin{\mathcal H}$, then $G$ is 3-factor-critical.
\end{theorem}
\noindent{\bf Proof.} By Lemma \ref{lem1}, diam$(G)=2$ or $3$. If
diam$(G)=3$,  then, by Lemma \ref{lem6}, the assertion follows.
We may now assume that diam$(G)=2$.

Suppose to the contrary that $G$ is not 3-factor-critical.
Then there exists a subset $S\subseteq V(G)$
such that $|S|=3$ and $G-S$ contains exactly two odd components by Lemma \ref{lem7}.
Let $S$, $C_{1}$ and $C_{2}$ be defined as before.
By Lemma \ref{lem8} (3), without loss of generality, suppose that $A_{1}\cap
A_{2}\neq \emptyset$. Further, it follows from Lemma \ref{lem9} that
both $\cap_{i=1}^{3}A_{i}=\emptyset$ and $\cap_{i=1}^{3}B_{i}=\emptyset$.
Thus there exist $i$ and $j$ such that $B_{i}\cap B_{j}\neq\emptyset$ by Lemma \ref{lem8} (4), where $1\leq i\neq j\leq 3$.
We next consider the following two subcases.

\noindent{\bf Case 1.} $B_{1}\cap B_{2}\neq \emptyset$.

Take $a\in A_{1}\cap A_{2}$, $b\in B_{1}\cap B_{2}$, and consider
the graph $G+ab$. If $|D_{ab}\cap \{a,b\}|=2$, then $|D_{ab}\cap S|=1$ by
Lemma \ref{lem5}. Since $\cap_{i=1}^{3}A_{i}=\emptyset$ and
$\cap_{i=1}^{3}B_{i}=\emptyset$, $D_{ab}\cap S=\{s_{1}\}$ or
$\{s_{2}\}$. This implies that $D_{ab}$ is a DDS of $G$, a
contradiction. Hence $|D_{ab}\cap \{a,b\}|=1$. By the symmetry of
structure of $G$, without loss of generality, we may assume that $a\in
D_{ab}$ and $b\notin D_{ab}$.

\noindent{\bf Case 1.1.} $|D_{ab}\cap S|=2$.

Clearly, we have $D_{ab}\cap S=\{s_{1},s_{3}\}$ or $\{s_{2},s_{3}\}$. Because the both cases of $D_{ab}\cap S$ can be discussed similarly, thus we may assume that
$D_{ab}\cap S=\{s_{1},s_{3}\}$.
Then each vertex of
$V(C_{2})-\{b\}$ is adjacent to $s_{1}$ and $s_{3}$. Hence
$(V(C_{2})-\{b\})\subseteq B_{1}\cap B_{3}$, and so
$V(C_{2})=B_{1}$. Take $b_{1}\in V(C_{2})-\{b\}$ and consider
the graph $G+ab_{1}$. Clearly, $s_{2}b_{1}\notin E(G)$ as
$\cap_{i=1}^{3}B_{i}=\emptyset$. If
$|D_{ab_{1}}\cap\{a,b_{1}\}|=2$, then $|D_{ab_{1}}\cap S|=1$ by
Lemma \ref{lem5}. Furthermore, note that $D_{ab_{1}}\cap S=\{s_{2}\}$ or $\{s_{3}\}$,
because $s_{1}\notin D_{ab_{1}}$.  If $s_{2}\in
D_{ab_{1}}$, then $(V(C_{2})-\{b,b_{1}\})\subseteq
\cap_{i=1}^{3}B_{i}\neq\emptyset$, a contradiction. So $s_{3}\in
D_{ab_{1}}$. But, $D_{ab_{1}}\nsucc\{b\}$ as $ab\notin E(G)$ and
$bs_{3}\notin E(G)$, a contradiction. Therefore,
$|D_{ab_{1}}\cap\{a,b_{1}\}|=1$.

Suppose $a\in D_{ab_{1}}$ and $b_{1}\notin D_{ab_{1}}$. If
$|D_{ab_{1}}\cap S|=2$, then $s_{2}\in D_{ab_{1}}$. However,
$D_{ab_{1}}\nsucc V(C_{2})-\{b,b_{1}\}$ as
$\cap_{i=1}^{3}B_{i}=\emptyset$. So $|D_{ab_{1}}\cap S|=1$. By
Lemma \ref{lem5}, in order to doubly dominate $V(C_{2})-\{b_{1}\}$,
 we have $|D_{ab_{1}}\cap (V(C_{2})-\{b_{1}\})|=1$. Note that $s_{3}\notin D_{ab_{1}}$,
 since $as_{3}\notin
E(G)$.  Recall that
$\cap_{i=1}^{3}B_{i}=\emptyset$. If $s_{2}\in D_{ab_{1}}$, then
$D_{ab_{1}}\nsucc V(C_{2})-\{b,b_{1}\}$. Hence $s_{1}\in
D_{ab_{1}}$. It immediately follows from Lemma \ref{lem8} (1) that $G[B_{1}]=G[V(C_{2})]$ is
complete. This implies that $D_{ab_{1}}$ is also a DDS of $G$, a
contradiction. Therefore, $a\notin D_{ab_{1}}$ and $b_{1}\in
D_{ab_{1}}$. Suppose $|D_{ab_{1}}\cap S|=2$. Then $D_{ab_{1}}\cap
S=\{s_{1},s_{3}\}$ or $\{s_{2},s_{3}\}$. Note that
$s_{2}b_{1}\notin E(G)$. If $D_{ab_{1}}\cap S=\{s_{1},s_{3}\}$,
then $s_{2}s_{1},s_{2}s_{3}\in E(G)$ to doubly dominate $s_{2}$.
This means that $S$ is a DDS of $G$, a contradiction. So
$D_{ab_{1}}\cap S=\{s_{2},s_{3}\}$. To doubly dominate $s_{2}$ and
$s_{1}$, $s_{2}s_{3}\in E(G)$ and $s_{1}$ is adjacent to at least
one of $s_{2}$ and $s_{3}$, respectively. This implies that $S$ is
a DDS of $G$, a contradiction. Hence $|D_{ab_{1}}\cap S|=1$. To
doubly dominate $V(C_{1})-\{a\}$, we see that $|D_{ab_{1}}\cap
(V(C_{1})-\{a\})|=1$ by Lemma \ref{lem5}. Since
$\cap_{i=1}^{3}B_{i}=\emptyset$, $bs_{3}\notin E(G)$. Then
$s_{3}\notin D_{ab_{1}}$. Recall that $s_{2}b_{1}\notin E(G)$.
Then $s_{2}\notin D_{ab_{1}}$. So $s_{1}\in D_{ab_{1}}$. Thus
$V(C_{1})=A_{1}$. By Lemma \ref{lem8} (1), $G[V(C_{1})]$ is complete, which
implies that $D_{ab_{1}}$ is a DDS of $G$, a contradiction. Hence
Case 1.1 can not occur.

\noindent{\bf Case 1.2.} $|D_{ab}\cap S|=1$.

Since $\cap_{i=1}^{3}
A_{i}=\emptyset$, $as_{3}\notin E(G)$. Then
$D_{ab}\cap S=\{s_{1}\}$ or $\{s_{2}\}$. Suppose $D_{ab}\cap
S=\{s_{1}\}$. To doubly dominate $V(C_{2})-\{b\}$, we have
$|D_{ab}\cap (V(C_{2})-\{b\})|=1$ by Lemma \ref{lem5}. Then
each vertex of $V(C_{2})$ is adjacent to $s_{1}$. Thus
$V(C_{2})=B_{1}$. By Lemma \ref{lem8} (1), $G[V(C_{2})]$ is complete. This
means that $D_{ab}$ is a DDS of $G$, a contradiction. So
$D_{ab}\cap S=\{s_{2}\}$. By a similar argument, one reaches the
same contradiction. Therefore, Case 1.2 is impossible.

\noindent{\bf Case 2.} $B_{1}\cap B_{3}\neq \emptyset$ or
$B_{2}\cap B_{3}\neq \emptyset$.

Suppose first that $B_{1}\cap B_{3}\neq \emptyset$. Choose $u\in
A_{1}\cap A_{2}$ and $v\in B_{1}\cap B_{3}$. Now consider $G+uv$.
We distinguish the following two subcases.

\noindent{\bf Case 2.1.} $|D_{uv}\cap \{u,v\}|=2$.

Then $|D_{uv}\cap S|=1$ by Lemma \ref{lem5}. Clearly, $s_{1}\notin
D_{uv}$. Thus $D_{uv}\cap S=\{s_{2}\}$ or $s_{3}$. First suppose
$D_{uv}\cap S=\{s_{2}\}$. Then each vertex of $V(C_{1})$ is
adjacent to $s_{2}$ while each vertex of $V(C_{2})-\{v\}$ is
adjacent to $s_{2}$ and $v$. By Lemma \ref{lem8} (1), $G[V(C_{1})]=G[A_{2}]$ is complete.
 Choose $v_{1}\in V(C_{2})-\{v\}$ and
consider $G+uv_{1}$. If $|D_{uv_{1}}\cap \{u,v_{1}\}|=2$, then
$|D_{uv_{1}}\cap S|=1$. Obviously, $s_{2}\notin D_{uv_{1}}$. If
$s_{1}\in D_{uv_{1}}$, then $(V(C_{2})-\{v,v_{1}\})\subseteq
B_{1}\cap B_{2}$. Thus $B_{1}\cap B_{2}\neq\emptyset$. By a
similar argument that used in the proof of Case 1, one reaches
the same contradictions. So $s_{1}\notin D_{uv_{1}}$ and $s_{3}\in
D_{uv_{1}}$. Then each vertex in $(V(C_{1})-\{u\})\cup V(C_{2})$
is adjacent to $s_{3}$. Since $us_{3}\notin E(G)$, $s_{2}s_{3}\in
E(G)$ to doubly dominate $s_{3}$ in $G+uv$. Note that
$\cap_{i=1}^{3}B_{i}=\emptyset$. So $s_{1}v_{1}\notin E(G)$ and
$s_{1}s_{3}\in E(G)$ to doubly dominae $s_{1}$ in $G+uv_{1}$. This
means that $S$ is a DDS of $G$, a contradiction. Hence
$|D_{uv_{1}}\cap \{u,v_{1}\}|=1$.

Suppose $u\in D_{uv_{1}}$ and $v_{1}\notin D_{uv_{1}}$. If
$|D_{uv_{1}}\cap S|=2$, then $s_{2}\notin D_{uv_{1}}$ as
$s_{2}v\notin E(G)$ and $uv\notin E(G)$. Thus
$D_{uv_{1}}=\{u,s_{1},s_{3}\}$ and each vertex of
$V(C_{2})-\{v,v_{1}\}$ is adjacent to both $s_{1}$ and $s_{3}$. Then
$(V(C_{2})-\{v,v_{1}\})\subseteq\cap_{i=1}^{3}B_{i}\neq\emptyset$,
a contradiction. So $|D_{uv_{1}}\cap S|=1$. To doubly dominate
$V(C_{2})-\{v_{1}\}$, we have $|D_{uv_{1}}\cap (V(C_{2})-\{v_{1}\})|=1$.
Recall that $s_{2}v\notin E(G)$ and $uv\notin E(G)$. Then
$s_{2}\notin D_{uv_{1}}$. Further $s_{3}\notin D_{uv_{1}}$ because
$us_{3}\notin E(G)$. Thus $s_{1}\in D_{uv_{1}}$. Then each vertex
of $V(C_{2})-\{v,v_{1}\}$ is adjacent to $s_{1}$, and so
$B_{1}\cap B_{2}\neq\emptyset$. By applying an argument similar to
that presented in the proof of Case 1, we can always reach a
contradiction. Therefore, $u\notin D_{uv_{1}}$ and $v_{1}\in
D_{uv_{1}}$. If $|D_{uv_{1}}\cap S|=2$, then
$D_{uv_{1}}=\{v_{1},s_{1},s_{3}\}$ or $\{v_{1},s_{2},s_{3}\}$.
Suppose $D_{uv_{1}}=\{v_{1},s_{1},s_{3}\}$. Then each vertex of
$V(C_{1})-\{u\}$ is adjacent to $s_{1}$ and $s_{3}$, which implies
that $\cap_{i=1}^{3}A_{i}\neq\emptyset$, a contradiction. Hence
$D_{uv_{1}}=\{v_{1},s_{2},s_{3}\}$. If $s_{1}$ is adjacent to a
vertex in $V(C_{2})-\{v\}$, then $B_{1}\cap B_{2}\neq\emptyset$.
Using a similar argument as in the proof of Case 1, one reaches
a contradiction. So $s_{1}$ is not adjacent to any vertex
in $V(C_{2})-\{v\}$. Similarly, $s_{3}$ is not adjacent to
any vertex in $V(C_{2})-\{v\}$. Thus $G-\{s_{2},v\}$ is not
connected, which contradicts the assumption that $G$ is
3-connected. Hence $|D_{uv_{1}}\cap S|=1$. To doubly dominate
$V(C_{1})-\{u\}$, we have $|D_{uv_{1}}\cap (V(C_{1})-\{u\})|=1$. Since
$G[V(C_{1})]$ is complete, $D_{uv_{1}}\cap
\{s_{1},s_{2}\}=\emptyset$ for otherwise $D_{uv_{1}}$ is a DDS of
$G$. So $s_{3}\in D_{uv_{1}}$. Then $A_{2}\cap A_{3}\neq\emptyset$
and $B_{2}\cap B_{3}\neq\emptyset$. By a similar argument that
used in the proof of Case 1, we can obtain a
contradiction. Hence $D_{uv}\cap S=\{s_{2}\}$ is impossible.
Similarly, $D_{uv}\cap S=\{s_{3}\}$ is also impossible. Therefore,
Case 2.1 can not occur.

\noindent{\bf Case 2.2.} $|D_{uv}\cap \{u,v\}|=1$.

\noindent{\bf Case 2.2.1.} $u\in D_{uv}$ and $v\notin
D_{uv}$.

If $|D_{uv}\cap S|=2$, then $D_{uv}=\{u,s_{1},s_{2}\}$ or
$\{u,s_{2},s_{3}\}$. Suppose $D_{uv}=\{u,s_{1},s_{2}\}$. Then each
vertex of $V(C_{2})-\{v\}$ is adjacent to both $s_{1}$ and $s_{2}$.
Thus $B_{1}\cap B_{2}\neq\emptyset$. By applying an argument
similar to that presented in the proof of Case 1, we can obtain a contradiction.
So $D_{uv}=\{u,s_{2},s_{3}\}$. Then each
vertex of $V(C_{2})-\{v\}$ is adjacent to both $s_{2}$ and $s_{3}$.
Thus $(V(C_{2})-\{v\})\subseteq B_{2}\cap B_{3}$ and
$G[V(C_{2})]=G[B_{3}]$ is complete by Lemma \ref{lem8} (1). Since $us_{3}\notin
E(G)$, $s_{2}s_{3}\in E(G)$. Further, $s_{1}$ is adjacent to at
least one of $s_{2}$ and $s_{3}$ because
$D_{uv}=\{u,s_{2},s_{3}\}\succ_{\times2}\{s_{1}\}$.

Take $v_{2}\in V(C_{2})-\{v\}$ and consider $G+uv_{2}$. If
$|D_{uv_{2}}\cap \{u,v_{2}\}|=2$, then $|D_{uv_{2}}\cap S|=1$.
Note that $uv,s_{2}v\notin E(G)$. So $s_{2}\notin D_{uv_{2}}$. If
$s_{1}\in D_{uv_{2}}$, then each vertex of $V(C_{2})-\{v,v_{2}\}$
is adjacent to $s_{1}$. Thus $\cap_{i=1}^{3}B_{i}\neq\emptyset$, a
contradiction. Hence $s_{1}\notin D_{uv_{2}}$ and $s_{3}\in
D_{uv_{2}}$. Then each vertex of $V(C_{1})-\{u\}$ is adjacent to
both $u$ and $s_{3}$. If there exists a vertex $u^{*}\in
V(C_{1})-\{u\}$ such that $u^{*}$ is adjacent to $s_{1}$ or
$s_{2}$, then $A_{1}\cap A_{3}\neq\emptyset$ or $A_{2}\cap
A_{3}\neq\emptyset$. Clearly, $B_{1}\cap B_{3}\neq\emptyset$ and
$B_{2}\cap B_{3}\neq\emptyset$. By a similar argument that used in
the proof of Case 1, we can obtain the same contradictions.
Hence $s_{1}$ and $s_{2}$ are not adjacent to any vertex in
$V(C_{1})-\{u\}$. Thus $G-\{u,s_{3}\}$ is not connected,
contradicting the fact that $G$ is 3-connected. Hence
$|D_{uv_{2}}\cap \{u,v_{2}\}|=1$.

Suppose $u\in D_{uv_{2}}$ and $v_{2}\notin D_{uv_{2}}$. If
$|D_{uv_{2}}\cap S|=2$, then $D_{uv_{2}}=\{u,s_{1},s_{2}\}$ or
$\{u,s_{1},s_{3}\}$. Since $uv,s_{2}v\notin E(G)$,
$D_{uv_{2}}=\{u,s_{1},s_{2}\}$ is impossible. So
$D_{uv_{2}}=\{u,s_{1},s_{3}\}$. Then each vertex of
$V(C_{2})-\{v,v_{2}\}$ is adjacent to $s_{1}$. Recall that
$(V(C_{2})-\{v\})\subseteq B_{2}\cap B_{3}$. Thus
$\cap_{i=1}^{3}B_{i}\neq\emptyset$, a contradiction. Hence
$|D_{uv_{2}}\cap S|=1$. Since $D_{uv_{2}}$ is a
$\gamma_{\times2}(G+uv_{2})$-set, $|D_{uv_{2}}\cap
(V(C_{2})-\{v_{2}\})|=1$. If $s_{2}\in D_{uv_{2}}$, then $vs_{2}\in E(G)$.
This derives that $\cap_{i=1}^{3}B_{i}\neq\emptyset$, a contradiction. Thus
$s_{2}\notin D_{uv_{2}}$. Further, $s_{3}\notin D_{uv_{2}}$ as
$us_{3}\notin E(G)$. So $s_{1}\in D_{uv_{2}}$. Then each vertex of
$V(C_{2})-\{v,v_{2}\}$ is adjacent to $s_{1}$, which implies that
$\cap_{i=1}^{3}B_{i}\neq\emptyset$, a contradiction. Therefore,
$u\notin D_{uv_{2}}$ and $v_{2}\in D_{uv_{2}}$. If
$|D_{uv_{2}}\cap S|=2$, then $D_{uv_{2}}=\{v_{2},s_{1},s_{3}\}$ or
$\{v_{2},s_{2},s_{3}\}$. Thus each vertex of $V(C_{1})-\{u\}$ is
adjacent to both $s_{1}$ and $s_{3}$ or both $s_{2}$ and $s_{3}$,
respectively. This means that $D_{uv_{2}}$ is a DDS of $G$, a
contradiction. Hence $|D_{uv_{2}}\cap S|=1$. To doubly dominate
$V(C_{1})-\{u\}$ in $G+uv_{2}$, we see that $|D_{uv_{2}}\cap
(V(C_{1})-\{u\})|=1$, say $u^{*}\in D_{uv_{2}}\cap
(V(C_{1})-\{u\})$. It is easy to see that $s_{2}\notin D_{uv_{2}}$
as $s_{2}v\notin E(G)$ and $u^{*}v\notin E(G)$. Since
$\cap_{i=1}^{3}B_{i}=\emptyset$, $s_{1}v_{2}\notin E(G)$, and so
$s_{1}\notin D_{uv_{2}}$. Therefore, $s_{3}\in D_{uv_{2}}$. Then
$s_{1}u^{*}\in E(G)$ and $s_{3}u^{*}\in E(G)$. Thus $A_{1}\cap
A_{3}\neq\emptyset$. Recall that $v\in B_{1}\cap B_{3}$. By a similar argument that used in the proof
of Case 1, we can get the same contradictions. Hence
$|D_{uv}\cap S|=2$ is impossible.

Thus $|D_{uv}\cap S|=1$. To doubly dominate $V(C_{2})-\{v\}$ in
$G+uv$, $|D_{uv}\cap (V(C_{2})-\{v\})|=1$, say $v^{*}\in
D_{uv}\cap (V(C_{2})-\{v\})$. If $s_{1}\in D_{uv}$, then
$V(C_{2})=B_{1}$. By Lemma \ref{lem8} (1), $G[V(C_{2})]$ is complete. This
implies that $D_{uv}$ is a DDS of $G$, a contradiction. So
$s_{1}\notin D_{uv}$. Since $us_{3}\notin E(G)$, $s_{3}\notin
D_{uv}$. Hence $s_{2}\in D_{uv}$. Then each vertex of
$V(C_{1})\cup (V(C_{2})-\{v\})$ is adjacent to $s_{2}$. Since
$D_{uv}\succ_{\times2}\{s_{3}\}$, $s_{3}v^{*}\in E(G)$. Thus $B_{2}\cap
B_{3}\neq\emptyset$. Further, $s_{3}$ is adjacent to a vertex of
$V(C_{1})-\{u\}$ because $A_{3}\neq\emptyset$. So $A_{2}\cap
A_{3}\neq\emptyset$. Using an argument similar to that presented
in the proof of Case 1, one reaches the same contradictions.
Therefore, Case 2.2.1 is impossible.

\noindent{\bf Case 2.2.2.} $u\notin D_{uv}$ and $v\in
D_{uv}$.

By applying a similar argument that used in the proof of Case
2.2.1, we can obtain the same contradictions. So Case 2.2 can
not occur. Hence $B_{1}\cap B_{3}=\emptyset$. Similarly, we
can show that $B_{2}\cap B_{3}=\emptyset$. Thus Case 2 can
not occur. This completes the proof of Theorem \ref{the1}.  ~$\Box$

\begin{figure}[ht]\label{fig1}
\setlength{\unitlength}{.037in}
\begin{picture}(110,75)
\put(32,9){\begin{picture}(110, 75)

\multiput(25,5)(60,0){2}{\circle*{1.8}}
\multiput(5,25)(40,0){2}{\circle*{1.8}}
\multiput(65,25)(40,0){2}{\circle*{1.8}}
\multiput(55,-5)(0,60){2}{\circle*{1.8}}
\put(85,25){\circle*{1.8}}

\put(65,25){\line(20,0){20}}\put(85,25){\line(20,0){20}}\put(5,25){\line(40,0){40}}
\put(25,5){\line(-1,1){20.2}}\put(25,5){\line(1,1){20.1}}
\put(85,5){\line(-1,1){20.2}}\put(85,5){\line(1,1){20.1}}
\put(55,-5){\line(-1,3){10}}\put(55,-5){\line(1,3){10}}
\put(55,-5){\line(-5,3){50}}\put(55,-5){\line(5,3){50}}
\put(55,55){\line(-5,-3){50}}\put(55,55){\line(5,-3){50}}
\put(55,55){\line(-3,-5){30}}\put(55,55){\line(3,-5){30}}
\put(55,55){\line(-1,-3){10.1}}\put(55,55){\line(1,-3){10.1}}
\put(55,55){\line(1,-1){29.6}}

\put(55,-5){\line(3,1){29.5}}\put(25,5){\line(3,1){59.5}}
\put(5,25){\line(4,-1){79.5}}\put(45,25){\line(2,-1){39.8}}

\put(55,-9){$z_2$}\put(55,57){$z_1$}\put(0,23){$x_{1}$}\put(25,2){$x_{3}$}
\put(47,24.5){$x_{2}$}\put(60,24){$y_{1}$}\put(107,24){$y_{2}$}
\put(85,2){$y_{3}$}\put(85,22){$y_{4}$}

\end{picture}}
\end{picture}
\caption{The graph $H_{6,3}$}
\end{figure}
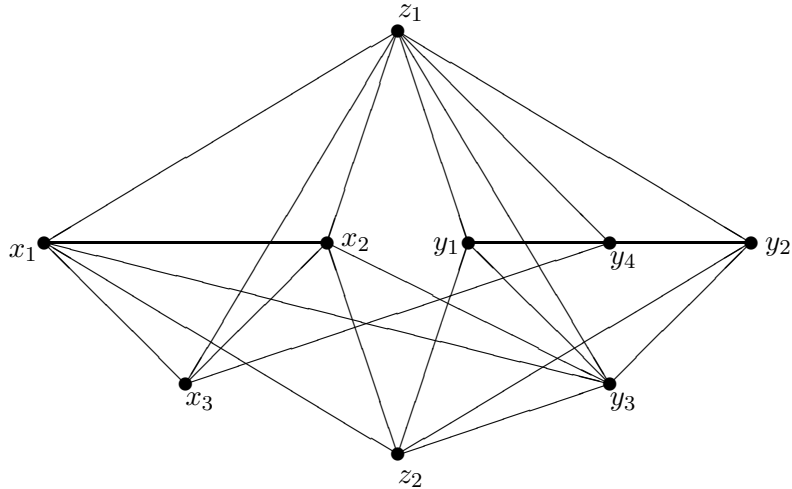

\noindent{\bf Remark.} The hypotheses on the
connectivity  and the
 minimum degree bounds in Theorem \ref{the1} are sharp. Indeed,
  the following results was proved by Favaron in \cite{f}: for all
 $k\geq 0$, every $k$-factor-critical graph of order $>k$ is
 $k$-(vertex)-connected and for all $k\geq 1$, every $k$-factor-critical
 graph of order $>k$ is $(k+1)$-(edge)-connected (and hence has
 minimum degree at least $k+1$). Next, to illustrate that the assumption on
 claw-free is necessary, we construct a graph $H_{6,t}$ as
 follows. For odd integer $t\geq 3$, let
 $F_1=K_t$, $F_2=C_4$ and $F_3=\overline{K}_2$. Moreover, let
 $V(F_1)=\{x_{1}, x_{2},\ldots
 x_{t}\}$, $V(F_2)=\{y_{1},y_{2},y_{3},y_{4}\}$ and $V(F_3)=\{z_1,z_2\}$.
 Let $H_{6,t}$ be the graph obtained from the union $F_1\cup F_2\cup F_3$ by
 adding edges $x_ty_4$ and $y_{3}x_{i}$ for $1\leq i\leq t-1$,
 joining $z_1$ to each vertex of $V(F_1)\cup V(F_2)$ and  joining $z_2$
 to each vertex of $(V(F_1)\cup V(F_2))-\{y_{4},x_{t}\}$.
 It is easy to see that $H_{6,t}$ is a 4-connected
 4-$\gamma_{\times2}$-critical graph of  order $t+6$ with minimum degree 4, but it
 contains a claw. Obviously, $H_{6,t}-\{u,v,y_{4}\}$ has no
 perfect matching. Hence $H_{6,t}$ is not 3-factor-critical. The graph $H_{6,3}$ is
 shown in Figure 1.

\end{document}